\documentclass[12pt]{article}
\usepackage{amssymb,amsfonts}
\usepackage{amscd}
\usepackage{amsmath}
\usepackage{graphicx}

\newtheorem{theorem}{Theorem}

\newtheorem{corollary}[theorem]{Corollary}

\newtheorem{example}[theorem]{Example}

\newtheorem{proposition}[theorem]{Proposition}

\newenvironment{proof}[1][Proof]{\noindent\textbf{#1.} }{\ \rule{0.5em}{0.5em}}

\begin{document}

\title{Cuntz-Pimnser Algebras,\\ Completely Positive Maps\\
and\\
Morita Equivalence}
\author{Alberto E. Marrero\thanks{Supported by grants from the National Science Foundation and the Sloan Foundation and by a GAANN Fellowship}\ and Paul S. Muhly\thanks{Supported by a grant from the National Science Foundation.} \medskip \\
%EndAName
Department of Mathematics\\
The University of Iowa\\
Iowa City, IA 52242\medskip\\
amarrero@math.uiowa.edu\\
pmuhly@math.uiowa.edu\\
}

\maketitle

\begin{abstract}
Let $P$ be a completely positive map on $M_n(\mathbb{C})$ and let $E_P$ be the associated \emph{GNS}-$C^*$-correspondence. We prove a result that implies, in particular, that the Cuntz-Pimsner algebra of $E_P$, $\mathcal{O}(E_P)$, is strongly Morita equivalent to the Cuntz algebra $\mathcal{O}_{d(P)}$, where $d(P)$ is the index of $P$.
\end{abstract}

By definition $E_{P}:=M_{n}\left( \mathbb{C}\right) \otimes _{P}M_{n}\left( \mathbb{C}%
\right) $ is the algebraic tensor product of $M_{n}\left( \mathbb{C}\right) $
with $M_{n}\left( \mathbb{C}\right) $ endowed with the $M_{n}\left( \mathbb{%
C}\right) $-valued inner product: $\langle a_{1}\otimes b_{1},a_{2}\otimes
b_{2}\rangle =b_{1}^{\ast }P\left( a_{1}^{\ast }a_{2}\right) b_{2}$, $%
a_{i}\otimes b_{i}\in E_{P}$, $i=1,2$, and obvious left and right actions of 
$M_{n}\left( \mathbb{C}\right) $ by matrix multiplication. Thus $E_P$ is a $C^{\ast }$-correspondence over $M_{n}(\mathbb{C})$ in the sense of \cite{MS98}.  In the literature, $E_{P}$ is often referred to as the \emph{GNS}-$C^*$-correspondence determined by $P$. 

The notion of index for $P$, $d(P)$, has a somewhat complicated definition \cite[Definition~4.9]{MS03}.  However, thanks to the special situation in which we are working, an alternative formulation is possible. To state it, recall that thanks to Stinespring's thereorm, we may write $P$ in the form
\begin{equation}\label{Stinespring}
P(a)=\sum\limits_{r=1}^{N}t_{r}at_{r}^{\ast },a\in M_{n}(\mathbb{C})
\end{equation}%
where $t_{1},t_{2},...t_{N}$ are in $M_{n}(\mathbb{C)}$. By Proposition 4.11 of \cite{MS03} the complex dimension of
the linear span of the $t_{i}$'s depends only on $P$ and is the index $d\left(
P\right) $ of $P$ defined in \cite[Definition~4.9]{MS03}.

The principal result of this note is

\begin{theorem}\label{ME_Correspondence}
The correspondence $E_P$ is strongly Morita equivalent, in the sense of \cite{MS00}, to the Hilbert space $\mathbb{C}^{d(P)}$.
\end{theorem}

\begin{proof} 
Since $M_{n}(\mathbb{C)}$ is strongly Morita equivalent to $%
\mathbb{C}$ via the column space $C_{n}(\mathbb{C})$, the
correspondence $E_{P}$ is strongly Morita equivalent to the correspondence
over $\mathbb{C}$, $F_{P}: =R_{n}(\mathbb{C})\otimes _{M_{n}(\mathbb{C}%
)}E_{P}\otimes _{M_{n}(\mathbb{C})}C_{n}(\mathbb{C})$, where $R_{n}\left( 
\mathbb{C}\right) $ is the $n$-dimensional row space over $\mathbb{C}$.  (See \cite{BMP00} for discussions of column and row spaces.) Thus $F_{P}$ is a finite dimensional Hilbert space over $\mathbb{C}$. We need to calculate its dimension. Lemma 4.10 of \cite{MS03} allows us to write $d(P)$ as $\dim _{\mathbb{C}}(e_{11}M_{n}(\mathbb{C})\otimes _{P}M_{n}(\mathbb{C}%
)e_{11})$,%
where $\{e_{ij}\}$ are the matrix units in $M_n(\mathbb{C})$.  Since the tensor products in the definition of $F_P$ as $R_{n}(\mathbb{C})\otimes _{M_{n}(\mathbb{C}%
)}E_{P}\otimes _{M_{n}(\mathbb{C})}C_{n}(\mathbb{C})$ are balanced over $M_{n}(\mathbb{C})$, the space $F_P$ is spanned by the elements $\varepsilon_{ij}:=e_{j}\otimes
I\otimes_{P}I\otimes e_{i}^{\ast}$ where $ e_{j}=\left( 0, \ldots ,1, \ldots ,0\right) $, with $1$ in the $j^{th}$ slot, and $i,j=1,\ldots ,n$. On the other hand, the space $(e_{11}M_{n}(\mathbb{C})\otimes _{P}M_{n}(\mathbb{C}%
)e_{11})$ is spanned by the matrices $e_{1i}{\otimes}_P e_{j1}$, $i,j=1,\ldots ,n$.  Consequently, to show that the two spaces, $F_P$ and $(e_{11}M_{n}(\mathbb{C})\otimes _{P}M_{n}(\mathbb{C}%
)e_{11})$, are isomorphic and so have the same dimension, it suffices to show that the inner products $\langle e_{1i}{\otimes}_P e_{j1}, e_{1k}{\otimes}_P e_{l1} \rangle$ and $\langle \varepsilon_{ij}, \varepsilon_{kl} \rangle$ are the same.  To do this, we write the matrices $t_r$ in the representation of $P$ in equation (\ref{Stinespring}) as $t_{r}=\left[ t\left( r;i,j\right) \right]_{i,j=1}^n$, 
$r=1, \ldots , N$.  Then, in $F_P$, we have
\begin{eqnarray*}
\langle \varepsilon_{ij}, \varepsilon_{kl} \rangle & = & \langle e_{i} \otimes I {\otimes}_P I \otimes e^*_{j}, e_{k} \otimes I {\otimes}_P I \otimes e^*_{l} \rangle \\ & = &%
e_j P(e^*_{i}e_{k})e^*_{l}  \\ & = & \sum_{r = 1}^{N} t(r; i,j) \overline{t(r; k,l)}. 
\end{eqnarray*}
On the other hand, by \cite[Lemma~4.10]{MS03}, the inner product in $(e_{11}M_{n}(\mathbb{C})\otimes _{P}M_{n}(\mathbb{C}%
)e_{11})$ is given by the formula $\left\langle a_{1}\otimes b_{1},a_{2}\otimes b_{2}\right\rangle =$ $%
tr\left( b_{1}^{\ast }P(a_{1}^{\ast }a_{2})b_{2}\right) $. So, we find that 
\begin{eqnarray*}
\langle e_{1i} \otimes e_{j1}, e_{1k} \otimes e_{l1} \rangle & = & tr(e^*_{j1} P(e^*_{1i} e_{1k}) e_{l1}) \\ & = & \sum_{r = 1}^{N} t(r; i,j) \overline{t(r; k,l)},
\end{eqnarray*}%
also.
\end{proof}

\begin{corollary}
Let $P$ be a completely positive map on $M_n(\mathbb{C})$ and let $E: = E_P$ be the associated $C^*$-correspondence over $M_n(\mathbb{C})$. Then
\begin{enumerate}
\item The tensor algebra over $E_P$ in the sense of \cite{MS98}, $\mathcal{T}_+(E_P)$ is strongly Morita equivalent, in the sense of \cite{BMP00} to Popescu's noncommutative disc algebra \cite{gP91}  $\mathcal{A}_{d(P)}$ of the size $d(P)$.
\item The Toeplitz $C^*$-algebra of $E_P$ is strongly Morita equivalent to the Toeplitz extension of Cuntz algebra $\mathcal{O}_{d(P)}$.
\item The Cuntz-Pimsner algebra $\mathcal{O}(E_P)$ is strongly Morita equivalent to the Cuntz algebra $\mathcal{O}_{d(P)}$. 
\end{enumerate}
\end{corollary}
\begin{proof}
The first two assertions are immediate from Theorem \ref{ME_Correspondence} and \cite[Theorem~3.2]{MS00}, once it is noted that the tensor algebra of $\mathbb{C}^{d(P)}$ is $\mathcal{A}_{d(P)}$. (See \cite{MS98}.) The third assertion follows from Theorem \ref{ME_Correspondence} and \cite[Theorem~3.5]{MS00}.
\end{proof}
\begin{example}
Let $P:M_{n}(\mathbb{C})\rightarrow M_{n}(\mathbb{C})$\ be the completely
positive map that takes a matrix to its diagonal part, i.e., suppose that $P(a) = \sum_{i} e_{ii}ae_{ii}$, $a \in M_{n}(\mathbb{C})$.  Evidently, the linear span of the $e_{ii}$ is $n$-dimensional and so we conclude that $\mathcal{\mathcal{O}}(E_{P})$ is strongly Morita
equivalent to $\mathcal{\mathcal{O}}_{n}$.
\end{example}
 
The argument used in the proof of Theorem \ref{ME_Correspondence} allows us to show what happens when one iterates the \emph{GNS}-construction with completely positive maps. Suppose we have two completely positives maps, $P_{1}$ and $P_{2}$,  acting 
on $M_{n}(\mathbb{C})$. Then we obtain a $C^*$-correspondence similar to the correspondence of the type $E_P$ we have been studying. It is the algebraic tensor product $M_{n}(%
\mathbb{C})\otimes_{P_{1}}M_{n}(C)\otimes_{P_{2}}M_{n}(\mathbb{C})$, with the obvious left and right multiplications by elements from $M_{n}(\mathbb{C})$ and $M_{n}(\mathbb{C})$-valued inner product given by the formula $\langle a_{1}\otimes b_{1}\otimes
c_{1},a_{2}\otimes b_{2}\otimes
c_{2}\rangle=c_{1}^{\ast}P_{2}(b_{1}^{\ast}P_{1}(a_{1}^{%
\ast}a_{2})b_{2})c_{2}$. We denote this correspondence by $E_{P_{1}P_{2}}$.

\begin{proposition}
Let $P_{1}$ and $P_{2}$ be two completely positive maps on $M_{n}(\mathbb{C%
})$. \ Then \ $E_{P_{1}P_{2}}$ is strongly
Morita equivalent to $\mathbb{C}^{d(P_1)d(P_2)}$.
\end{proposition}

\begin{proof}
Since \ $M_{n}(\mathbb{C}%
) $ is strongly Morita equivalent to $\mathbb{C}$, via the column space $%
C_{n}(\mathbb{C})$ ,the correspondence $F_{P_{1}P_{2}}:=R_{n}(\mathbb{C}%
)\otimes_{M_{n}(\mathbb{C})}E_{P_{1}P_{2}}\otimes _{M_{n}(\mathbb{C})}C_{n}(%
\mathbb{C})$ is strongly Morita equivalent to $E_{P_{1}P_{2}}$. Likewise $F_{P_{i}}:=R_{n}(\mathbb{C})\otimes_{M_{n}(\mathbb{C}%
)}E_{P_{i}}\otimes_{M_{n}(\mathbb{C})}C_{n}(\mathbb{C})$ is strongly Morita equivalent to $E_{P_{i}}$, $i = 1, 2$.  We show that $F_{P_{1}P_{2}}$ is naturally isomorphic to $F_{P_1} \otimes F_{P_2}$.

Observe that $F_{P_{1}P_{2}}$ is generated  by $\varepsilon_{iklj}:=e_{j}%
\otimes I\otimes e_{kl}\otimes I\otimes e_{i}^{\ast}$ \. Further, the inner product in  $F_{P_{1}P_{2}}$ is given by the formula
\begin{align*}
\langle\xi_{1}\otimes I\otimes a\otimes I\otimes\eta_{1},\xi_{2}\otimes
I\otimes b\otimes I\otimes\eta_{2}\rangle_{F_{P_{1}P_{2}}} & =
\eta_{1}^{\ast}P_{2}(a^{\ast}P_{1}(\xi_{1}^{\ast}\xi_{2})b)\eta_{2}.
\end{align*}
From this, a simple computation along the lines that we made in the proof of Theorem \ref{ME_Correspondence}, where, recall, $\varepsilon_{ij} = e_j \otimes I \otimes_P I \otimes e_k^*$, shows that $\left\langle\varepsilon_{iklj},%
\varepsilon_{mphq}\right\rangle _{F_{P_{1}P_{2}}}=\left\langle
\varepsilon_{kj},\varepsilon_{pq}\right\rangle _{F_{P_{1}}}\left\langle
\varepsilon_{il},\varepsilon_{mh}\right\rangle _{F_{P_{2}}} = \langle \varepsilon_{kj} \otimes \varepsilon_{il}, \varepsilon_{pq} \otimes \varepsilon_{mk}\rangle$. This, in turn, implies that $F_{P_{1}P_{2}}$ is isomorphic to $F_{P_{1}} \otimes F_{P_{2}}$.
\end{proof}

\bigskip

\noindent \textbf{Acknowledgement}
\emph{The first author's contribution to this note was made while he was a visiting scholar at the School of Mathematical and Physical Sciences at the University of Newcastle, Australia.  He would like to thank the School for its hospitality and, especially, he would like to thank Iain Raeburn his continuous encouragement and generous advice.}


\begin{thebibliography}{9}
\bibitem{BMP00} D. Blecher, P. Muhly and V. Paulsen, \textit{Categories of
operator modules (Morita equivalence and projective modules),} Memoirs of
the Amer. Math. Soc., Vol. 143, \# 681, Providence, 2000.

\bibitem{cL95} E.Christopher Lance, \textit{Hilbert C*-Modules,} London
Mathematical Society, Lecture Notes Series 210.

\bibitem{MS98} P. Muhly and B. Solel, \textit{Tensor Algebras over
C*-Correspondences: Representations, Dilations, and C*-Envelopes,} J.
Functional Anal. \textbf{158 }(1998), 389--457.

\bibitem{MS00} P. Muhly and B. Solel, \textit{On the Morita
Equivalence of Tensor Algebras,} Proc. London Math Soc. \textbf{81} (2000),
113--168.

\bibitem{MS03} P. Muhly and B. Solel, \textit{The Curvature and
index of completely positive maps, }Proc. London Math. Soc. \textbf{87}
(2003), 748--778.

\bibitem{gP91}  G. Popescu, Von Neumann inequality for $(B(H)^{n})_{1}$, Math.Scand. \textbf{68} (1991), 292--304. 


\end{thebibliography}
\end{document}